\documentclass[12pt]{article}
\usepackage{amsmath}
\usepackage{amssymb,amsgen,amsbsy,amsopn,amsfonts,graphicx,theorem}
\usepackage[dvips]{color}
\usepackage{multicol}
\usepackage{mathrsfs}

\numberwithin{equation}{section}

  \font\tenmsb=msbm10
     \font\sevenmsb=msbm7
     \font\fivemsb=msbm5
     \catcode`\@=11
     \ifx\amstexloaded@\relax\catcode`\@=\active
     \endinput\else\let\amstexloaded@\relax\fi
     \def\spaces@{\space\space\space\space\space}
     \def\spaces@@{\spaces@\spaces@\spaces@\spaces@\spaces@}
     \def\space@.{\futurelet\space@\relax}
     \space@. %
     \def\Err@#1{\errhelp\defaulthelp@\errmessage{AmS-TeX error: #1}}
     \def\relaxnext@{\let\next\relax}
     \def\accentfam@{7}
     \def\noaccents@{\def\accentfam@{0}}
     \def\Cal{\relaxnext@\ifmmode\let\next\Cal@\else
     \def\next{\Err@{Use \string\Cal\space only in math mode}}\fi\next}
     \def\Cal@#1{{\Cal@@{#1}}}
     \def\Cal@@#1{\noaccents@\fam\tw@#1}
     \def\Bbb{\relaxnext@\ifmmode\let\next\Bbb@\else
     \def\next{\Err@{Use \string\Bbb\space only in math mode}}\fi\next}
     \def\Bbb@#1{{\Bbb@@{#1}}}
     \def\Bbb@@#1{\noaccents@\fam\msbfam#1}
     \newfam\msbfam
     \textfont\msbfam=\tenmsb
     \scriptfont\msbfam=\sevenmsb
     \scriptscriptfont\msbfam=\fivemsb

\newtheorem{theorem}{Theorem}[section]

\newtheorem{condition}{Condition}[section]

\newtheorem{definition}{Definition}[section]

\newtheorem{lemma}{Lemma}[section]

\newtheorem{remark}{Remark}[section]

\def\v{\varepsilon}

\begin{document}
\title{Vanishing Viscosity Limit for Isentropic Navier-Stokes Equations with Density-dependent  Viscosity}
\author{ {Feimin Huang$^a$, Ronghua Pan$^b$, Tianyi Wang$^{a}$,
 Yong Wang$^{a,}$\footnote{Corresponding author.  \newline \indent
Email addresses: fhuang@amt.ac.cn(F. Huang),
panrh@math.gatech.edu(R.Pan), wangtianyi@amss.ac.cn(T. Wang),
yongwang@amss.ac.cn(Y. Wang),  zhaixy@amt.ac.cn(X. Zhai)}, Xiaoyun
Zhai$^{a}$ }
\\
\ \\
   {\small \it $^a$Institute of Applied Mathematics,  Academy of Mathematics and Systems
   Science}\\
{\small \it Chinese Academy of Sciences, Beijing, 100190, China} \\
{\small \it $^b$School of  Mathematics, Georgia Institute of
Technology, Atlanta, GA 30332}}
\date{ }

\maketitle

%\vskip 0.1cm \arraycolsep1.5pt

\begin{abstract}
In this paper, we study the vanishing viscosity limit of
one-dimensional isentropic compressible Navier-Stokes equations with density-dependent viscosity, to the
isentropic compressible Euler equations.
Based on several new uniform estimates to the viscous systems, in addition to the framework recently established by G. Chen
and M. Perepelitsa \cite{Chen6}, we justify that the finite energy solution of
the isentropic compressible Euler equations for a large class of  initial data can be obtained as the inviscid limit of the
compressible Navier-Stokes equations even when the viscosity depends on the density.

 \


\

AMS: 35L50, 35L60, 35L65, 76R50
\end{abstract}

\section{Introduction}

When the fluid density experiences large scale dropping,
especially when vacuum is concerned, the motion of isentropic
compressible viscous fluids
is modeled by the following compressible
Navier-Stokes equations with the density-dependent viscosity, in
the Eulerian coordinates,
\begin{eqnarray}\label{1.1}
\begin{cases}
\rho^{\varepsilon}_{t}+({\rho^{\varepsilon} u^{\varepsilon}})_x=0,\\
({\rho^{\varepsilon} u^{\varepsilon}})_t+(\rho^\varepsilon
(u^\varepsilon)^2+p(\rho^\varepsilon))_x
=\varepsilon((\rho^{\varepsilon})^{\alpha}u^{\varepsilon}_{x})_x,
\end{cases}
\end{eqnarray}
where $\rho^{\varepsilon}$ and
$u^{\varepsilon}$ denote the density and the velocity of the fluid,
respectively. $m^\varepsilon=\rho^\varepsilon u^\varepsilon$
represents the momentum.  $p=p(\rho)$ is pressure function of the
density. In this paper, we consider the polytropic perfect gas, i.e.
$$p\left({\rho}\right)=\kappa\rho^{\gamma},$$
where $\gamma > 1$ is the adiabatic exponent, and
the constant $\kappa$ is chosen
as $\kappa=\frac{\left({\gamma-1}\right)^2}{4\gamma}$ up to a scaling. 
While $\varepsilon>0$ is adpated to
the system as the controlling parameter on the amplitude of
viscosity, for which we assume
$\varepsilon\in\left({0,\varepsilon_0}\right]$ for some fixed
$\varepsilon_0>0$ without loss of generality; $\alpha\ge 0$ is
a constant which models the dependence
of viscosity on density.

When $\alpha=0$, (1.1) reduces to the classical
compressible Navier-Stokes equation, called {\bf CNS}. The case of $\alpha>0$ occurs for
non-uniform gases \cite{Chapman}, and (1.1) can be formally
derived by Chapman-Enskog expansion from the Boltzmann equation for (at least) hard sphere model and cut-off
inverse power force model. A formal derivation can be found in
\cite{MW}. It is also interesting to note that when $\alpha=1$
and $\gamma=2$, (1.1) recovers the ``viscous Saint-Venant" system
for shallow water without bottom friction \cite{shallowwater}, see also
\cite{Courant-Friedrichs}.
In this paper, we will
focus on (1.1) with positive $\alpha$, for which we call it
{\bf $\alpha$-CNS},
distinguishing from the case of $\alpha=0$, which is called {\bf CNS}.
On the other hand, the studies in \cite{Hoff-Serre} and \cite{liu-xin-yang} indicate the failure of CNS at vacuum and the validity of $\alpha$-CNS
at least at the level of local well-posedness theory. We therefore devote
our efforts to this model in current paper.

We now consider the Cauchy problem of (1.1) when the far fields of
the fluid are away from vacuum. Namely, we shall study the $\alpha$-CNS
(1.1) with the following initial data
\begin{eqnarray}\label{1.2}
\rho^{\varepsilon}(0,x)=\rho^{\varepsilon}_0(x)>0 ,\quad
u^{\varepsilon}(0,x)=u^{\varepsilon}_0(x),
\end{eqnarray}
such that 
$$\lim_{x\rightarrow
\pm\infty}(\rho^{\varepsilon}_0(x),u^{\varepsilon}_0(x))=(\rho^\pm,u^\pm),
\ with \ \rho^\pm > 0.$$

In the past decades, the study of the mathematical theory on (1.1)--(1.2) has attracted a lot attention. Many interesting results were established for the
local and global existence of both classical and weak solutions, we refer the
readers to some of them such as, \cite{BD1}, \cite{BDG}, \cite{Fang-Zhang1}, \cite{Fang-Zhang2},
\cite{Hoff}, \cite{Song Jiang}, \cite{JX}, \cite{JXZ}, \cite{Kanel},
\cite{LLX}, \cite{liu-xin-yang}, \cite{MW}, \cite{MV1}, \cite{MV2} and
\cite{Zhu changjiang}. It is equally
interesting to study the inviscid limit for (1.1)--(1.2) as $\varepsilon\to 0$
toward the following one-dimensional isentropic Euler equations
\begin{eqnarray}\label{1.3}
\begin{cases}
\rho_{t}+(\rho u)_x=0,\\
(\rho u)_t+(\rho u^2+p(\rho))_x=0.
\end{cases}
\end{eqnarray}
It is a general belief that the physcial weak solution of (1.3) can
be obtained in such a process, see \cite{A. Bressan}, where a vanishing artificial viscosity
limit for general hyperbolic system with small BV data is proved. 
This problem is closely related to the
existence of weak solutions to (1.3) through a limitting process of
physical approximation. In this paper, we will address
this problem and study the vanishing viscosity limit for
(1.1)--(1.2).

In BV framework, when the initial data is away from vacuum, the existence of global BV solution to (1.3) was established by \cite{Nishida-Smoller} for $\gamma>1$ and by \cite{Nishida} for
$\gamma=1$ using Glimm's method. (1.3) shows singular behavior 
when vacuum occurs which
causes difficulties to mathematical analysis. It is still a major open
problem on how to perform BV estimate when the solution may contain vacuum states. Instead, the $L^\infty$ framework is successfully achieved using the
theory of compensated compactness \cite{Murat}, \cite{L. Tartar}. The existence
of $L^\infty$ weak entropy solution of (1.3) was established by
\cite{R.J.DiPerna1} for
$\gamma=1+\frac{2}{2n+1}, n\geq2$; by \cite{Chen5} for $\gamma\in(1,\frac53]$;
by \cite{Lions P.-L.1} and \cite{Lions P.-L.2} for $\gamma>\frac53$; and finally
by \cite{Feimin Huang-Zheng Wang} for $\gamma=1$. Recently,
\cite{Ph. LeFloch} further constructed
the finite-energy solutions to the isentropic Euler equations with
finite-energy initial data. We remark that these results are achieved
through the vanishing artificial viscosity.

The problem of vanishing physical viscosity limit is more
subtle and the progress has been less satisfactory, and the problem of vanishing viscosity limit of Navier-Stokes equations to
Euler equations has been open for long time, though some interesting results
are proved when restrictive initial data is assigned, see \cite{Hoff-Liu}
and \cite{Xin}. Recently, G.
Chen and M. Perepelitsa \cite{Chen6} proved that the solutions of Navier-Stokes
($\alpha=0$),
whose viscosity is independent of density, converge to the finite
energy solution of Euler equations as viscosity vanishes. This is a
major breakthrough in this aspect.

Inspired by \cite{Chen6}, we
study the problem of vanishing viscosity for the $\alpha$-CNS (1.1)--(1.2)
in this paper with positive $\alpha$. It is clear that for any fixed positive
$\varepsilon$, the visocity coefficient with positive $\alpha$ experiences
degeneracy near vacuum states. An obvious obstacle is the dissipation term 
in the energy identity contains only the weighted norm of velocity gradient which degenerates at vacuum. Such a singular behavior causes the major
difficulty in the analysis and introduced the different behavior of solutions
compared with CNS where $\alpha=0$. The analysis exibits quite different flavor and requires very different ingredients. 
Fortunately, by a deep observation, we obtained several key uniform
estimates. Based on these uniform estimates and the framework of
\cite{Chen6}, we are able to show that, when viscosity parameter $\v$ tends to
zero, the solutions of $\alpha$-CNS
\eqref{1.1}-(1.2) converge to the
finite-energy solution of Euler equations for general initial data.

We now prepare to state our main result.

A pair of functions $(\eta(\rho,u),q(\rho,u))$, or
 $(\eta(\rho,m),q(\rho,m))$ for
$m=\rho u$, is called an entropy-entropy flux pair of system (1.3), if the following holds
$$[\eta(\rho,u)]_t+[q(\rho,u)]_x=0,$$
for any smooth solutions of (1.3). Furthermore, $\eta(\rho,m)$
is called a weak entropy if
$$\eta(0, u)=0, \ \mbox{for any fixed}\ u.$$

An entropy $\eta(\rho, m)$ is convex if the Hessian
$\nabla^2\eta(\rho,m)$ is nonnegative definite in the region under
consideration.

From \cite{Lions P.-L.2}, it is well known that any  week entropy
$\left({\eta,q}\right)$ can be represented by
\begin{equation}\label{1.4}
\begin{cases}
&\displaystyle \eta^\psi(\rho,\rho u)=\eta^\psi(\rho,m)=\int_{\mathbb{R}}{\chi\left({\rho;s-u}\right)\psi\left({s}\right)ds},\\
&\displaystyle q^\psi(\rho,\rho u)=q^\psi(\rho,m)=\int_{\mathbb{R}}\left({\theta
s+\left({1-\theta}\right)u}\right)\chi\left({\rho;s-u}\right)\psi\left({s}\right)ds.
\end{cases}
\end{equation}
where the kernel is
$\chi(\rho;s-u)=[\rho^{2\theta}-({s-u})^2]^{\lambda}_+$, 
$\lambda=\frac{3-\gamma}{2(\gamma-1)}>-\frac12$, and $\theta=\frac{\gamma-1}{2}$. For instance, when $\psi(s)=\frac{1}{2}s^2$, the entropy pair is the
mechanical energy and the associated flux
\begin{eqnarray}
\eta^{\ast}(\rho,m)=\frac{m^2}{2\rho}+e(\rho),\quad
q^{\ast}(\rho,m)=\frac{m^3}{2\rho^2}+m e'(\rho),
\end{eqnarray}
where $e(\rho)=\frac{\kappa}{\gamma-1}\rho^\gamma$ represents the gas
internal energy in physics.

\

 Let
$\left({\bar\rho\left({x}\right),\bar{u}\left({x}\right)}\right)$ be
a pair of smooth monotone functions satisfying
$\left({\bar\rho\left({x}\right),\bar{u}\left({x}\right)}\right)=\left({\rho^\pm,u^\pm}\right)$,
when $\pm x \geq  L_0$ for some large $L_0 > 0 $. The total
mechanical energy for (\ref{1.1}) in $\mathbb{R}$ with respect to
the pair of reference function
$\left({\bar\rho\left({x}\right),\bar{u}\left({x}\right)}\right)$ is
\begin{eqnarray}\label{1.6}
E\left[{\rho,u}\right]\left({t}\right)=\int_{\mathbb{R}}{\left({\eta^{\ast}\left({\rho,m}\right)
-\eta^{\ast}\left({\bar\rho,\bar{m}}\right)-
\nabla\eta^{\ast}(\bar\rho,\bar{m})\cdot\left({\rho-\bar\rho,m-\bar{m}}\right)
}\right)dx},
\end{eqnarray} where $\bar{m}=\bar\rho\bar{u}$. After some
calculations, we obtain that
\begin{eqnarray}\label{1.7}
E\left[{\rho,u}\right]\left({t}\right)
=\int_{\mathbb{R}}\left({\frac{1}{2}\rho\left({t,x}\right)\left|{u(t,x)-\bar{u}(x)}\right|^2
+e^{\ast}\left({\rho(t,x),\bar\rho(x)}\right)}\right)dx
\end{eqnarray}
where
$e^{\ast}\left({\rho,\bar\rho}\right)=e(\rho)-e(\bar\rho)-e'\left({\bar\rho}\right)\left({\rho-\bar\rho}\right)\geq
0$.

\begin{definition}\label{Def1.1}
Let $\left({\rho_0,u_0}\right)$ be given initial data with
finite-energy with respect to the end states $\left({\rho^\pm
,u^\pm}\right)$ at infinity, and $E\left[{\rho_0 ,u_0}\right]\le E_0
<\infty$. A pair of measurable functions $\left({\rho ,u}\right):
\mathbb{R}^2_+ \rightarrow \mathbb{R}^2_+$ is called a finite-energy
entropy solution of the Cauchy problem (1.3) if the following
holds:

(i) The total energy in bounded in time: There is a bounded function
$C\left({E,t}\right)$, defined on $\mathbb{R}^+\times\mathbb{R}^+$
and continuous in $t$ for each $E\in \mathbb{R}^+$, such that, for
a.e. $t>0$, $$E\left[{\rho, u}\right]\left({t}\right)\le
C\left({E_0, t}\right);$$

(ii) The entropy inequality:
$$\eta^{\psi}\left({\rho,u}\right)_t+q^{\psi}\left({\rho,u}\right)_x\le0,$$
is satisfied in the sense of distributions for all test functions
$\psi\left({s}\right)\in \left\{{\pm 1,\pm s,s^2}\right\}$;

(iii) The initial data $\left({\rho_0,u_0}\right)$ are attained in
the sense of distributions.
\end{definition}

We now state our main conditions on the initial data (1.2), which is 
motivated from \cite{Chen6}.

\begin{condition} Let $\left({\bar\rho\left({x}\right),\bar{u}\left({x}\right)}\right)$
be some pair of smooth monotone functions satisfying
$\left({\bar\rho\left({x}\right), 
 \bar{u}\left({x}\right)}\right)=\left({\rho^{\pm}\left({x}\right), u^{\pm}\left({x}\right)}\right)$
 when $\pm x \geq L_0$ for some large $L_0>0$. For positive constants
 $E_0$, $E_1$ and $M_0$ independent of $\varepsilon$, and $c^{\varepsilon}_0>0$
the initial functions
$\left({\rho^{\varepsilon}_0,u^{\varepsilon}_0}\right)$ are smooth
satisfying the following properties

(i) $\rho^{\varepsilon}_0\geq c^{\varepsilon}_0>0,\quad
\displaystyle\int_{\mathbb{R}}\rho^{\varepsilon}_0\left({x}\right)\left|{u^{\varepsilon}_0\left({x}\right)-\bar{u}\left({x}\right)}\right|dx\leq
M_0<\infty$ ;

(ii) The total mechanical energy with respect to
$\left({\bar\rho,\bar{u}}\right)$ is finite:
$$\int_{\mathbb{R}}\left({\frac{1}{2}}\rho^{\varepsilon}_0
\left|{u^{\varepsilon}_0\left({x}\right)-\bar{u}\left({x}\right)}\right|^2+e^{\ast}\left({\rho^{\varepsilon}_0\left({x}\right),\bar\rho\left({x}\right)}\right)\right)dx\leq
E_0<\infty;$$

(iii)
$\displaystyle\varepsilon^2\int_{\mathbb{R}}
\frac{\left|{\rho^{\varepsilon}_{0,x}\left({x}\right)}\right|^2}{\rho^{\varepsilon}_{0}
\left({x}\right)^{3-2\alpha}}dx\leq
E_1<\infty$;

(iv)$\left({\rho^{\varepsilon}_0\left({x}\right),\rho^{\varepsilon}_0\left({x}\right)u^{\varepsilon}_0
\left({x}\right)}\right)\rightarrow
\left({\rho_0\left({x}\right),\rho_0\left({x}\right)u_0\left({x}\right)}
\right)$
in the sense of distributions as $\varepsilon \rightarrow 0$, with
$\rho_0\left({x}\right) \geq 0\quad a.e.$. 

\end{condition}

Our main results are stated in the following Theorem.

\begin{theorem}\label{th1.1}
Assume $\frac23\leq\alpha\leq\gamma$, $\gamma>1$. Let $\left({\rho^{\varepsilon},u^{\varepsilon}}\right)$, $m^{\varepsilon}=\rho^{\varepsilon}u^{\varepsilon}$ be the solution of the Cauchy problem (\ref{1.1})-(1.2) with initial data $\left({\rho^{\varepsilon}_0,u^{\varepsilon}_0}\right)$ which satisfies 
Condition 1.1 for each fixed $\varepsilon > 0$. Then,
when $\varepsilon\rightarrow 0$, there exists a subsequence of
$\left({\rho^{\varepsilon},m^{\varepsilon}}\right)$  that converges
almost everywhere to a finite-energy entropy solution
$\left({\rho,m}\right)$ to the Cauchy problem (\ref{1.3}) with
initial data $\left({\rho_0,\rho_0 u_0}\right)$ for the isentropic
Euler equations.
\end{theorem}

\begin{remark} Due to some technical difficulty, we can only prove the 
result for $\frac23\leq\alpha\leq\gamma$. In fact, for many physical gases
the Chapman-Enskog 
viscosity predicts that $\alpha\ge \frac{\gamma-1}{2}$, see \cite{Chapman}, 
and \cite{MW}. Our condition $\frac23\leq\alpha\leq\gamma$ is valid for many physical cases including the shallow water model, but it did not cover
the case of monoatomic gas where $\gamma=\frac53$ and $\alpha=\frac12$.
It is very interesting to prove the result for $\alpha\in[0,\frac23]$, which
will be addressed later.
\end{remark}
One important basis of our proof for Theorem \ref{th1.1} is the following
compactness theorem,  established in  \cite{Chen6}.

\begin{theorem}[Chen-Perepelitsa \cite{Chen6}]\label{th4.1}
Let $\psi\in C_0^2(\mathbb{R})$, $(\eta^{\psi},q^{\psi})$ be a weak
entropy pair generated by $\psi$. Assume that the sequences
$(\rho^{\varepsilon}(x,t),u^{\varepsilon}(x,t))$ defined on
$\mathbb{R}\times\mathbb{R}_+$ with
$m^{\varepsilon}=\rho^{\varepsilon}u^{\varepsilon}$, satisfies the
following conditions:

 (i).   For any $-\infty<a<b<\infty$ and all
$t>0$, it holds that
\begin{eqnarray}\label{3.201}
\int_{0}^{t}\int_{a}^{b}(\rho^\varepsilon)^{\gamma+1}dxd\tau\leq
C(t,a,b),
\end{eqnarray}
where $C(t)>0$ is  independent of $\varepsilon$.

(ii).For any compact set $K\subset\mathbb{R}$, it holds that
\begin{eqnarray}\label{3.301}
\int_{0}^{t}\int
_{K}(\rho^\varepsilon)^{\gamma+\theta}+\rho^\varepsilon|u^\varepsilon|^3dxd\tau\leq
C(t,K),
\end{eqnarray}
where $C=C(t,K)>0$ is independent of $\varepsilon$.

(iii). The sequence of entropy dissipation measures
\begin{eqnarray}\label{4.1-1}
\eta^{\psi}(\rho^{\varepsilon},m^{\varepsilon})
_t+q^{\psi}(\rho^{\varepsilon},m^{\varepsilon})_x \ \mbox{are
compact in $H_{loc}^{-1}(\mathbb{R}_+^2)$}.
\end{eqnarray}
Then there is a subsequence of 
$(\rho^{\varepsilon},m^{\varepsilon})$(still denoted as $(\rho^{\varepsilon},m^{\varepsilon}))$ and a pair of
measurable functions $(\rho, m)$ such that
\begin{eqnarray}\label{4.1-2}
(\rho^{\varepsilon},m^{\varepsilon})\rightarrow (\rho, m),\ \ a.e. \ as\ 
\varepsilon\rightarrow0.
\end{eqnarray}
\end{theorem}

In section 2 and section 3 below, we will verify conditions
\eqref{3.201}, \eqref{3.301} and \eqref{4.1-1}  to prove our main
theorem \ref{th1.1}.

The rest of this paper is arranged as follows. In section 2, we make
some new uniform estimates for the solutions of Navier-Stokes
equations (\ref{1.1}) which are independent of $\varepsilon$. These
estimates are essential to show the convergence of the
vanishing viscosity limit to the Euler equations. In section 3, 
using the estimate we obtained in section 2, we prove the $H_{loc}^{-1}(\mathbb{R}_+^2)-$ compactness  for the solutions of
(\ref{1.1}). In section 4, based on the framework in
\cite{Chen6}, we prove our main Theorem\ref{th1.1}.

%%%%%%%%%%%%%%%%%%%%%%%%%%%%%%%%%%%%%%%%%%%%%%%%%%%%%%%%%%%%%%%%%%%%%%%%%%%%%%%%%%%%%%%%%%%%%%%%%%%%%%
\section{Uniform Estimates for the Solutions of $\alpha$-CNS}\label{section3}

First, we assume that $(\rho^\varepsilon,u^\varepsilon)$ is the
global smooth solutions of Navier-Stokes equations
\eqref{1.1}--\eqref{1.2}, satifying
\begin{eqnarray}
\rho^{\varepsilon}(x,t)\geq c^{\varepsilon}(t), \ \mbox{for some}\
c^{\varepsilon}(t)>0
\end{eqnarray}
and
\begin{eqnarray}
\lim_{x\rightarrow \pm\infty}(\rho^{\varepsilon},
u^\varepsilon)(x,t)=(\rho^{\pm}, u^{\pm}).
\end{eqnarray}
For the existence of global smooth solutions, the reders are referred to
\cite{JX}, \cite {MW} and \cite{MV2}.
Based on the above preparation, we now make some new the uniform
estimates with respect to $\varepsilon$ for the solutions $(\rho^\varepsilon,u^\varepsilon)$ of
the $\alpha$-CNS \eqref{1.1}--\eqref{1.2}. 

For simplicity, throughout this section, we denote
$(\rho,u)=(\rho^\varepsilon,u^\varepsilon)$ without causing confusion
and $C>0$ denote the constant independent of $\varepsilon$.

\begin{lemma}(Energy Estimates)\label{lem3.1}
Suppose that  $0\leq \alpha \leq \gamma$, and $E\left[{\rho_0,u_0}\right]\leq
E_0<\infty$ for some $E_0>0$ independent of $\varepsilon$. It holds
that
\begin{eqnarray}\label{3.1}
\sup_{0 \leq \tau\leq
t}  E[\rho,u](\tau)+\varepsilon\int_{0}^{t}\int_{\mathbb{R}}\rho^{\alpha}u_{x}^2
dxd\tau\leq C(t),
\end{eqnarray}
where $C(t)$ depends on $E_0$, $t$, $\bar\rho$, and $\bar{u}$, but not on $\varepsilon$.
\end{lemma}

\noindent\textbf{Proof}. From the definition, we have
\begin{eqnarray}\label{3.2}
\frac{dE(t)}{dt}=\frac{d}{dt}\int_{\mathbb{R}}\eta^{\ast}(\rho,m)-\eta^{\ast}(\bar\rho,\bar
m)\ dx
-\int_{\mathbb{R}}\nabla\eta^{\ast}(\bar\rho,\bar{m})(\rho_{t},m_t)\ dx.
\end{eqnarray}

Since $(\eta^{\ast},q^{\ast})$ is an entropy pair, we have
\begin{eqnarray}\label{3.3}
\eta^{\ast}(\rho,m)_t+q^{\ast}(\rho,m)_x
-\varepsilon\eta_m^{\ast}(\rho,m)(\rho^\alpha u_x)_x=0.
\end{eqnarray}
Integrate \eqref{3.3} with respect to $x$ over $\mathbb{R}$, we
obtain
\begin{eqnarray}\label{3.4}
\frac{d}{dt}\int_{\mathbb{R}}\eta^{\ast}(\rho,m)-\eta^{\ast}(\bar\rho,\bar
m)\ dx +\varepsilon\int_{\mathbb{R}}\rho^\alpha
u_x^2\ dx=q^{\ast}(\rho^-,m^-)-q^{\ast}(\rho^+,m^+).
\end{eqnarray}

Since we have
\begin{eqnarray}\label{3.50}
e^{\ast}(\rho, \bar\rho)\geq\rho(\rho^\theta-\bar{\rho}^\theta)^2,\
\theta=\frac{\gamma-1}{2}.
\end{eqnarray}

Utilizing \eqref{3.50}, we obtain
\begin{eqnarray}\label{3.5}
\begin{aligned}
&|\int_{\mathbb{R}}\nabla\eta^{\ast}(\bar\rho,\bar{m})(\rho_{t},m_t)dx|\nonumber\\
&=|\int_{\mathbb{R}}\nabla\eta^{\ast}(\bar\rho,\bar{m})(m_x,(p(\rho)+\rho
u^2-\varepsilon\rho^\alpha u_x)_x)dx|\nonumber\\
&=|\int_{\mathbb{R}}\nabla\eta^{\ast}(\bar\rho,\bar{m})_x(m,p(\rho)+\rho
u^2-\varepsilon\rho^\alpha u_x)\ dx|\nonumber\\
&\leq \frac{\varepsilon}{4}\int_{\mathbb{R}}\rho^{\alpha}u_{x}^2
\ dx+C\int_{\mathbb{R}}\rho|u-\bar{u}|^2\ dx+C\int_{-L_0}^{L_0}(\rho+p(\rho)+\rho^\alpha)dx+C\nonumber\\
&\leq
C+CE+\frac{\varepsilon}{4}\int_{\mathbb{R}}\rho^{\alpha}u_{x}^2 dx,
\end{aligned}
\end{eqnarray}
where we have used
\begin{eqnarray}\label{3.6}
\int_{-L_0}^{L_0}(\rho+p(\rho)+\rho^\alpha)dx\leq CE, \ \ \mbox{for}
\ \ 0\leq\alpha\leq \gamma.
\end{eqnarray}

Substituting \eqref{3.4} and \eqref{3.5} into \eqref{3.2}, we obtain
\begin{eqnarray}\label{3.7}
\frac{dE(t)}{dt}+\frac{3\varepsilon}{4}\int_{\mathbb{R}}\rho^{\alpha}u_{x}^2\ dx\leq
C+CE,
\end{eqnarray}
Then Gronwall's inequality implies Lemma$\ref{lem3.1}$.

\begin{remark}
Since vacuum could occur in our solution, the inequality 
$$\int_{0}^{t}\int_
{\mathbb{R}}\rho^{\alpha}u_{x}^2\ dxd\tau\leq C(t)$$
in \eqref{3.1} is
much weaker than the corresponding one
$$\int_{0}^{t}\int_
{\mathbb{R}}u_{x}^2 dxd\tau\leq C(t).$$ 
in \cite{Chen6}. This will cause a great
difficulty to prove Lemma \ref{lem3.3} below, which is an essential
step to verify the condition i) of Theorem \ref{th4.1}, i.e.
\eqref{3.201}.
\end{remark}

We now derive some higher order estimates.

\begin{lemma}\label{lem3.2}
If $0<\alpha \leq \gamma$, and $(\rho_0(x), u_0(x))$ satisfies
$$\varepsilon^2\int_{\mathbb{R}}\frac{|\rho_{0x}(x)|^2}{\rho_{0}(x)^{3-2\alpha}}\ dx\leq
E_1<\infty,$$ 
for some $E_1$ independent of $\varepsilon$. Then, for any
$t>0$, it holds that
\begin{eqnarray}\label{3.8}
\varepsilon^2\int_{\mathbb{R}}\rho^{2\alpha-3}\rho_x^2\ dx
+\varepsilon\int_{0}^{t}\int_{\mathbb{R}}\rho^{\alpha+\gamma-3}\rho_{x}^2
\ dxd\tau\leq C(t),
\end{eqnarray}
where $C(t)$ depends on $E_0$, $E_1$, $t$, $\bar\rho$, $\bar{u}$, but
not on $\varepsilon$.
\end{lemma} 

\noindent\textbf{Proof}. Through \eqref{1.1}, we have
\begin{eqnarray}\label{3.9}
\begin{cases}
\rho_{xt}+\rho_{xx}u+2\rho_x u_x+\rho u_{xx}=0,\\
\rho u_t+\rho
uu_x+p(\rho)_x=\varepsilon\rho^{\alpha}u_{xx}
+\varepsilon\alpha\rho^{\alpha-1}u_x\rho_x.
\end{cases}
\end{eqnarray}

Multiplying $(\ref{3.9})_1$ with $\rho^{2\alpha-3}\rho_x$, after some
calculation, we obtain

\begin{eqnarray}\label{3.10}
\bigg(\frac{\rho^{2\alpha-3}\rho_x^2}{2}\bigg)_t
+\bigg(\frac{\rho^{2\alpha-3}u\rho_x^2}{2}\bigg)_x
+\alpha\rho^{2\alpha-3}\rho_x^2u_x+\rho^{2\alpha-2}\rho_xu_{xx}=0.
\end{eqnarray}

From $(\ref{3.9})_2\times\rho^{\alpha-2}\rho_x$, after some
calculation, we reach
\begin{eqnarray}\label{3.11}
\rho^{\alpha-1}\rho_xu_t
+\rho^{\alpha-1}u\rho_xu_x+\kappa\rho^{\alpha+\gamma-3}\rho_x^2
=\varepsilon\alpha\rho^{2\alpha-3}\rho_x^2u_x+\varepsilon\rho^{2\alpha-2}
\rho_xu_{xx}.
\end{eqnarray}

The combination $\varepsilon^2\eqref{3.10}+\varepsilon\eqref{3.11}$ gives
that
\begin{eqnarray}\label{3.12}
\begin{aligned}
&\bigg(\frac{\varepsilon^2\rho^{2\alpha-3}\rho_x^2}{2}\bigg)_t
+\bigg(\frac{\varepsilon^2\rho^{2\alpha-3}u\rho_x^2}{2}\bigg)_x
+\varepsilon\kappa\rho^{\alpha+\gamma-3}\rho_x^2\\
&\qquad\qquad+\bigg(\varepsilon\rho^{\alpha-1}\rho_xu\bigg)_t
-\bigg(\varepsilon\rho^{\alpha-1}\rho_tu\bigg)_x
=\varepsilon\rho^{\alpha}u_x^2.
\end{aligned}
\end{eqnarray}
Integrating \eqref{3.12} over $[0,t]\times\mathbb{R}$, we obtain
\begin{eqnarray}\label{3.13}
\begin{aligned}
&\varepsilon^2\int_{\mathbb{R}}\frac{\rho^{2\alpha-3}\rho_x^2}{2}\ dx
+\varepsilon\kappa\int_{0}^{t}\int_{\mathbb{R}}
\rho^{\alpha+\gamma-3}\rho_{x}^2\ 
dxd\tau\\
&=\varepsilon^2\int_{\mathbb{R}}\frac{\rho_0^{2\alpha-3}\rho_{0x}^2}{2}dx
+\varepsilon\int_{0}^{t}\int_{\mathbb{R}}\rho^{\alpha}u_x^2dxd\tau
-\varepsilon\int_{0}^{t}\int_{\mathbb{R}}
\bigg(\rho^{\alpha-1}\rho_xu\bigg)_{\tau}dxd\tau\\
&\leq
C(t)-\varepsilon\int_{0}^{t}\int_{\mathbb{R}}
\bigg(\rho^{\alpha-1}\rho_xu\bigg)_{\tau}dxd\tau.
\end{aligned}
\end{eqnarray}

Noticing that
\begin{eqnarray}\label{3.14}
\bigg(\rho^{\alpha-1}\rho_xu\bigg)_t
=\bigg(\rho^{\alpha-1}\rho_x(u-\bar{u})\bigg)_t
+\frac{1}{\alpha}\bigg((\rho^{\alpha})_x\bar{u}\bigg)_t,
\end{eqnarray}
Integrating \eqref{3.14} with respect to $x$ over $\mathbb{R}$, we
have
\begin{eqnarray}\label{3.15}
\begin{aligned}
&\varepsilon\int_{0}^{t}\int_{\mathbb{R}}
\bigg(\rho^{\alpha-1}\rho_xu\bigg)_{\tau}\ dxd\tau\\
&=\varepsilon\int_{\mathbb{R}}\rho^{\alpha-1}\rho_x(u-\bar{u})\ dx
-\varepsilon\int_{\mathbb{R}}\rho_0^{\alpha-1}\rho_{0x}(u_0-\bar{u})\ 
dx\\
&\qquad+\frac{\varepsilon}{\alpha}\int_{\mathbb{R}}(\rho^{\alpha})_x\bar{u}dx
+\frac{\varepsilon}{\alpha}\int_{\mathbb{R}}(\rho_0^{\alpha})_x\bar{u}\ dx\\
&\leq C(t),
\end{aligned}
\end{eqnarray}
where we have used the following estimates \eqref{3.16}-\eqref{3.18}
\begin{eqnarray}\label{3.16}
\begin{aligned}
&\varepsilon\int_{\mathbb{R}}\rho^{\alpha-1}\rho_x(u-\bar{u})\ dx\\
&\leq\frac{\varepsilon^2}{8}\int_{\mathbb{R}}\rho^{2\alpha-3}\rho_x^2dx
+C\int_{\mathbb{R}}\rho(u-\bar{u})^2\ dx\\
&
\leq\frac{\varepsilon^2}{8}\int_{\mathbb{R}}\rho^{2\alpha-3}\rho_x^2dx+C(t),
\end{aligned}
\end{eqnarray}
\begin{eqnarray}
&\displaystyle\varepsilon\int_{\mathbb{R}}\rho_0^{\alpha-1}\rho_{0x}(u_0-\bar{u})\ dx
\leq\frac{\varepsilon^2}{8}\int_{\mathbb{R}}\rho_0^{2\alpha-3}
\rho_{0x}^2dx+C(t),\label{3.17}\\
&\displaystyle \frac{\varepsilon}{\alpha}\int_{\mathbb{R}}(\rho^{\alpha})_x\bar{u}\ dx
=-\frac{\varepsilon}{\alpha}\int_{\mathbb{R}}\rho^{\alpha}\bar{u}_xdx
+\frac{\varepsilon}{\alpha}\Big((\rho^+)^{\alpha}\bar{u}^+ -
(\rho^-)^{\alpha}\bar{u}^-\Big)\nonumber\\
&\displaystyle\leq \frac{C\varepsilon}{\alpha}\int_{-L_0}^{L_0}\rho^{\alpha}dx+C(t)\nonumber\\
&\leq C(t),\label{3.18}
\end{eqnarray}
Substituting \eqref{3.15}--\eqref{3.18} into \eqref{3.13}, we obtain
Lemma$\ref{lem3.2}$.

The following higher order integrability estimate is crucial in compactness argument.

\begin{lemma}\label{lem3.3}
 If the conditions of
Lemma $\ref{lem3.1}$ hold and $0< \alpha \leq \gamma$, then
for any $-\infty<a<b<\infty$ and all $t>0$, it holds that
\begin{eqnarray}\label{3.20}
\int_{0}^{t}\int_{a}^{b}\rho^{\gamma+1}\ dxd\tau\leq C(t,a,b),
\end{eqnarray}
where $C(t)>0$ depends on $E_0$, $a$, $b$, $\gamma$, $t$, $\bar\rho$, 
 $\bar{u}$, but not on $\varepsilon$.
\end{lemma}

\noindent\textbf{Proof}. Choose $$w(x)\in C_0^{\infty}(\mathbb{R}),\ 
0\leq w(x) \leq 1, \ w(x)=1 \ \mbox{for}\ x\in [a,b], \ \mbox{and}\
 supp\{ w\}=(a-1,b+1)$$

By $\eqref{1.1}_2\times w$, we have
\begin{eqnarray}\label{3.21}
(p(\rho)w)_x=-(\rho u^2w)_x+(p(\rho)+\rho u^2)w_x-(\rho
u)_tw+\varepsilon(\rho^{\alpha}u_xw)_x-\varepsilon\rho^{\alpha}u_xw_x
\end{eqnarray}

Integrating \eqref{3.21} with respect to spatial variable over $(-\infty, x)$, we
obtain
\begin{eqnarray}\label{3.22}
\begin{aligned}
p(\rho)w&\displaystyle=-\rho
u^2w+\varepsilon\rho^{\alpha}u_xw-\bigg(\int_{-\infty}^{x}\rho
uw\ dy\bigg)_t\\
&\displaystyle+\int_{-\infty}^{x}[(\rho
u^2+p(\rho))w_x-\varepsilon\rho^{\alpha}u_xw_x]\ dy.
\end{aligned}
\end{eqnarray}

Multiplying \eqref{3.22} by $\rho w$, we have
\begin{equation}\label{3.23}
\begin{aligned}
\rho p(\rho)w^2&\displaystyle=-\rho^2
u^2w^2+\varepsilon\rho^{\alpha+1}u_x w^2-\bigg(\rho
w\int_{-\infty}^{x}\rho uw\ dy\bigg)_t\\
&\displaystyle\qquad-(\rho u)_x w\int_{-\infty}^{x}\rho uw\ dy+\rho
w\int_{-\infty}^{x}[(\rho
u^2+p(\rho))w_x-\varepsilon\rho^{\alpha}u_xw_x]\ dx\\
&\displaystyle=\varepsilon\rho^{\alpha+1}u_xw^2-\bigg(\rho
w\int_{-\infty}^{x}\rho uw\ dy\bigg)_t-\bigg(\rho u w\int_{-\infty}^{x}\rho uw\ dy\bigg)_x\\
&\displaystyle\qquad+\rho u w_x\int_{-\infty}^{x}\rho uw\ dy+\rho
w\int_{-\infty}^{x}[(\rho
u^2+p(\rho))w_x-\varepsilon\rho^{\alpha}u_xw_x]\ dy
\end{aligned}
\end{equation}

Integrating \eqref{3.23} over $(0,t)\times\mathbb{R}$, we have
\begin{equation}\label{3.24}
\begin{aligned}
&\displaystyle\int_{0}^{t}\int_{\mathbb{R}}\kappa\rho^{\gamma+1}w^2\ dxd\tau\\
&\displaystyle=\varepsilon\int_{0}^{t}\int_{\mathbb{R}}\rho^{\alpha+1}u_xw^2\ dxd\tau
-\int_{\mathbb{R}}\bigg(\rho w\int_{-\infty}^{x}\rho
uw\ dy\bigg)\ dx\\
&\displaystyle\quad+\int_{\mathbb{R}}\bigg(\rho_0 w\int_{-\infty}^{x}\rho_0
u_0wdy\bigg)dx+\int_{0}^{t}\int_{\mathbb{R}}\bigg(\rho u
w_x\int_{-\infty}^{x}\rho uw\ dy\bigg)dxd\tau\\
&\displaystyle\quad+\int_{0}^{t}\int_{\mathbb{R}}\bigg(\rho
w\int_{-\infty}^{x}[(\rho
u^2+p(\rho))w_x-\varepsilon\rho^{\alpha}u_xw_x]dy\bigg)\ dxd\tau.
\end{aligned}
\end{equation}

Let
\begin{eqnarray}\label{3.54}
A=\{x: \rho(x,t)\geq\hat\rho\},\ \mbox{where}\
\hat{\rho}=2\max\{\rho +,\rho -\},
\end{eqnarray}
then we have the following estimates by \eqref{3.1}
\begin{eqnarray}
|A|\leq \frac{C(t)}{e^\ast(2\hat\rho,\bar\rho)}=: d(t).
\end{eqnarray}

By \eqref{3.54}, we know that for any $(x,t)$ there  exists a point
$x_0=x_0(x,t)$ such that $|x-x_0|\leq d(t)$ and
$\rho(x_0,t)=\hat\rho$. Here we choose
$\beta=\alpha+\frac{\gamma-1}{2}>0$,
\begin{eqnarray}\label{3.55}
\begin{aligned}
&\displaystyle\sup_{x\in supp \{w\}}\varepsilon\rho^\beta(x,t)
 \leq \varepsilon\hat\rho^\beta+\sup_{x\in supp\{w\}\cap
 A}\varepsilon\rho^\beta(x,t)\\
&\displaystyle  \ \leq
2\varepsilon\hat\rho^\beta+\sup_{x\in supp\{ w\}\cap
 A}|\varepsilon\rho^\beta(x,t)-\varepsilon\rho^\beta(x_0,t)|\\
&\displaystyle \ \leq2\varepsilon\hat\rho^\beta
+\sup_{x\in supp\{w\}\cap A}
\int_{x_0-d(t)}^{x_0+d(t)}|\beta||\varepsilon\rho^{\beta-1}
\rho_x| \ dx\\
&\displaystyle \leq2\varepsilon\hat\rho^\beta
+\int_{a-1-2d(t)}^{b+1+2d(t)}|\beta||\varepsilon\rho^{\beta-1}\rho_x|\ dx\\
&\displaystyle\ \leq2\varepsilon\hat\rho^\beta
+\int_{a-1-2d(t)}^{b+1+2d(t)}|\beta|\rho^{2\beta-2\alpha+1}\ dx
+\int_{\mathbb{R}}\varepsilon^2\rho^{2\alpha-3}\rho_x^2\ dx\\
&\displaystyle \ \leq C(t)
+\int_{a-1-2d(t)}^{b+1+2d(t)}\rho^{\gamma}\ dx\\
&\displaystyle\ \leq C(t)
\end{aligned}
\end{eqnarray}

Using \eqref{3.55}, the first term on the right hand side of
\eqref{3.24} can be estimated as following
\begin{equation}\label{3.25}
\begin{aligned}
&\displaystyle\varepsilon\int_{0}^{t}\int_{\mathbb{R}}
\rho^{\alpha+1}u_xw^2dxd\tau\\
&\displaystyle\leq
\varepsilon\int_{0}^{t}\int_{\mathbb{R}}\rho^{\alpha+2}w^4\ dxd\tau
+\varepsilon\int_{0}^{t}\int_{\mathbb{R}}\rho^{\alpha}u_x^2\ dxd\tau\\
&\displaystyle\leq C(t)+\varepsilon\int_{0}^{t}\int_{\mathbb{R}}\rho^{\alpha+2}w^2\ dxd\tau\\
&\displaystyle\leq
\begin{cases}
\displaystyle 
 C(t)+\varepsilon\int_{0}^{t}\int_{\mathbb{R}}\rho^{\beta}w^2\ dxd\tau,
\
\mbox{if\ \ $2\alpha+2\leq\beta$ }\\
\displaystyle C(t)+\int_{0}^{t}\sup_{x\in supp\{w\}}
\varepsilon\rho^\beta(x,\tau)\int_{\mathbb{R}}\rho^{2\alpha+2-\beta}w^2\ 
 dxd\tau,\
\ \mbox{if\ \ $2\alpha+2>\beta$ }\\
\end{cases}\\
&\displaystyle\leq
C(t)+C(t)\int_{0}^{t}\int_{\mathbb{R}}\rho^{2\alpha+2-\beta}w^2\ dxd\tau\\
&\displaystyle\leq
C(t)+\delta\int_{0}^{t}\int_{\mathbb{R}}\rho^{\gamma+1}w^2\ dxd\tau.
\end{aligned}
\end{equation}
Here we have used the fact  $2\alpha+2-\beta<\gamma+1$ for $\gamma>1$.

By Lemma $\ref{lem3.1}$ and the H\"older inequality, we obtain
\begin{eqnarray}\label{3.26}
|\int_{-\infty}^{x}\rho u wdy|&\leq&\int_{supp \{w\}}|\rho u|
dy\nonumber\\
&=&\bigg(\int_{supp \{w\}}\rho dy\bigg)^\frac{1}{2}\bigg(\int_{supp
\{w\}}\rho u^2dy\bigg)^\frac{1}{2}\leq C(t).
\end{eqnarray}
Then it follows that
\begin{eqnarray}\label{3.27}
&\displaystyle|\int_{\mathbb{R}}\bigg(\rho w\int_{-\infty}^{x}\rho
uw\ dy\bigg)\ dx|+|\int_{\mathbb{R}}\bigg(\rho_0
w\int_{-\infty}^{x}\rho_0
u_0w\ dy\bigg)\ dx|\\
&\displaystyle+|\int_{0}^{t}\int_{\mathbb{R}}\bigg(\rho u
w_x\int_{-\infty}^{x}\rho uw\ dy\bigg)\ dxd\tau|\leq C(t).
\end{eqnarray}
Similarly, we have
\begin{eqnarray}\label{3.28}
|\int_{0}^{t}\int_{\mathbb{R}}\bigg(\rho w\int_{-\infty}^{x}(\rho
u^2+p(\rho))w_x\ dy\bigg)\ dxd\tau|\leq C(t),
\end{eqnarray}
and
\begin{equation}\label{3.29}
\begin{aligned}
&\displaystyle|\varepsilon\int_{0}^{t}\int_{\mathbb{R}}\bigg(\rho
w\int_{-\infty}^{x}\rho^{\alpha}u_xw_x\ dy\bigg)\ dxd\tau|\\
&\displaystyle\leq\varepsilon\int_{0}^{t}\int_{\mathbb{R}}\bigg(\rho
w\int_{\mathbb{R}}\rho^{\alpha}|u_x||w_x|\ dy\bigg)\ dxd\tau\\
&\displaystyle\leq \varepsilon\int_{0}^{t}\bigg(\int_{\mathbb{R}}\rho
w\ dx\bigg)\bigg(\int_{\mathbb{R}}\rho^{\alpha}u_x^2\ dy
+\int_{\mathbb{R}}\rho^{\alpha}w_x^2\ dy\bigg)\ d\tau\\
&\leq C(t)
\end{aligned}
\end{equation}

 Substituting \eqref{3.25}, $\eqref{3.27}-\eqref{3.29}$ into
 \eqref{3.24} and noticing the smallness of $\delta$, we proved
 Lemma$\ref{lem3.3}$.

\begin{lemma}\label{lem3.4}
Suppost that $\frac23\leq\alpha \leq \gamma$, and
$\left({\rho_0(x),u_0(x)}\right)$ satisfy the conditions in the
Lemmas \ref{lem3.1}, \ref{lem3.2}, \ref{lem3.3}. Furthermore, assume that
for some
$M_0>0$ independent of $\varepsilon$, the following 
\begin{eqnarray}\label{3.51}
\int_{\mathbb{R}}\rho_0(x)\left|{u_0(x)-\bar{u}(x)}\right|dx\leq
M_0<\infty,
\end{eqnarray} 
holds. Then for any compact set $K\subset\mathbb{R}$, it
holds that
\begin{eqnarray}\label{3.30}
\int_{0}^{t}\int_{K}\rho^{\gamma+\theta}+\rho|u|^3\ dxd\tau\leq
C(t,K),
\end{eqnarray}
where $C=C(E_0,E_1,M_0,t,K)>0$ is independent of $\varepsilon$.
\end{lemma}

\noindent\textbf{Proof}. First, we introduce a useful result about
the entropy pair, see \cite{Lions P.-L.2} for details. Taking
$\psi^{\star}(w)=\frac{1}{2}w|w|$, then
 there exists a positive constant $C>0$, depending only on
 $\gamma>1$, such that the corresponding entropy pair
$(\eta^{\star},q^{\star})=(\eta^{\psi^\star},q^{\psi^\star})$
satisfies
\begin{eqnarray}\label{3.52}
\begin{cases}
|\eta^{\star}(\rho,u)|\leq\left({\rho|u|^2+\rho^\gamma}\right),\\
q^{\star}(\rho,u)\geq
C^{-1}\left({\rho|u|^3+\rho^{\gamma+\theta}}\right),\\
|\eta^{\star}_m(\rho,u)|\leq C\left({|u|+\rho^{\theta}}\right),\\
|\eta^{\star}_{mm}(\rho,u)|\leq C \rho^{-1},
\end{cases}
 \mbox{for all
$\rho\geq0$ and $u\in\mathbb{R}$}.
\end{eqnarray}
If  $\eta^{\star}_m$ is  regarded as the function of  $(\rho, u)$,
we have
\begin{eqnarray}\label{3.53}
\begin{cases} |\eta^{\star}_{mu}(\rho,u)|\leq C,\\
|\eta^{\star}_{m\rho}(\rho,u)|\leq C\rho^{\theta-1},
\end{cases}
 \mbox{for all
$\rho\geq0$ and $u\in\mathbb{R}$}.
\end{eqnarray}

For this weak entropy pair $(\eta^{\star},q^{\star})$, we note that
$$\eta^{\star}(\rho,0)=\eta_{\rho}^{\star}(\rho,0)=0,
\ \
q^{\star}(\rho,0)=\frac{\theta}{2}\rho^{3\theta+1}
\int_{\mathbb{R}}|s|^3[1-s^2]_{+}^{\lambda}\ ds$$
and
$$\eta_{m}^{\star}(\rho,0)=\beta\rho^{\theta},\  \mbox{with}\  \beta:=\int_{\mathbb{R}}|s|[1-s^2]_{+}^{\lambda}\ ds .$$

Taylor expansion implies
\begin{eqnarray}\label{3.31}
\eta^{\star}(\rho,m)=\beta\rho^{\theta}m+r(\rho,m),
\end{eqnarray}
with
\begin{eqnarray}\label{3.32}
r(\rho,m)\leq C\rho u^2.
\end{eqnarray}
for some constant $C>0$. Now we introduce a new
entropy pair $(\hat{\eta},\hat{q})$,
$$\hat{\eta}(\rho,m)=\eta^{\star}(\rho,m-\rho u^-),\ \
\hat{q}(\rho,m)=q^{\star}(\rho,m-\rho
u^-)+u^{-}\eta^{\star}(\rho,m-\rho u^-),$$
which satisfies
\begin{eqnarray}\label{3.33}
\begin{cases}
\hat\eta(\rho,m)=\beta\rho^{\theta+1}(u-u^-)+r(\rho,\rho(u-u^-)),\\
r(\rho,\rho(u-u^-))\leq C\rho(u-u^-)^2.
\end{cases}
\end{eqnarray}

Integrating
$\eqref{1.1}_1\times\hat{\eta}_{\rho}+\eqref{1.1}_2\times\hat{\eta}_{m}$
 over $(0,t)\times(-\infty,x)$, we obtain
\begin{equation}\label{3.34}
\begin{aligned}
&\displaystyle\int_{-\infty}^{x}(\hat{\eta}(\rho,m)-\hat{\eta}(\rho_0,m_0))\ dy
+\int_{0}^{t}q^{\star}(\rho,\rho(u-\bar{u}))+u^{-}\eta^{\star}(\rho,\rho(u-\bar{u}))\ d\tau\\
&\displaystyle=tq^{\star}(\rho^-,0)+\varepsilon\int_{0}^{t}\hat{\eta}_m\rho^{\alpha}u_{x}\ d\tau
-\varepsilon\int_{0}^{t}\int_{-\infty}^{x}
(\hat{\eta}_{mu}\rho^{\alpha}u_{x}^2+ \hat{\eta}_{m\rho}\rho^{\alpha}
\rho_xu_{x})\ dyd\tau.
\end{aligned}
\end{equation}
Utilizing \eqref{3.53}, we have the following estimates
\begin{eqnarray}
&\displaystyle|\varepsilon\int_{0}^{t}\int_{-\infty}^{x}\hat{\eta}_{mu}
\rho^{\alpha}u_{x}^2dyd\tau|\leq
C
\varepsilon\int_{0}^{t}\int_{\mathbb{R}}\rho^{\alpha}u_{x}^2dyd\tau\leq
C(t),\label{3.35}\\
&\displaystyle|\varepsilon\int_{0}^{t}\int_{-\infty}^{x}
\hat{\eta}_{m\rho}\rho^{\alpha}\rho_xu_{x}dyd\tau|
\leq C\varepsilon\int_{0}^{t}\int_{\mathbb{R}}\rho^{\theta-1}\rho^{\alpha}|\rho_xu_{x}|\ dyd\tau\nonumber\\
&\displaystyle\leq
C\varepsilon\int_{0}^{t}\int_{\mathbb{R}}\rho^{\alpha}u_{x}^2dyd\tau
+C\varepsilon\int_{0}^{t}\int_{\mathbb{R}}\rho^{\alpha+\gamma-3}\rho_{x}^2\ dyd\tau\leq C(t)\label{3.36}.
\end{eqnarray}

Substituting \eqref{3.35} and \eqref{3.36} into \eqref{3.34}, then
integrating the result with respect to $x$ over $K$  and  using
\eqref{3.52}, we obtain
\begin{equation}\label{3.37}
\begin{aligned}
&\displaystyle\int_{0}^{t}\int_{K}\rho^{\alpha+\gamma}+\rho|u-u^{-}|^3\  
 dxd\tau\\
&\displaystyle\leq
C(t)+C\int_{0}^{t}\int_{K}|\eta^{\star}(\rho,\rho(u-\bar{u}))|\ dxd\tau
+C\varepsilon\int_{0}^{t}\int_{K}\rho^{\alpha}|u_x||u|\ dxd\tau\\
&\displaystyle
+C\varepsilon\int_{0}^{t}\int_{K}\rho^{\alpha+\theta}|u_x|\ dxd\tau\
+2\sup_{\tau\in[0,t]}\bigg|\int_{K}\bigg(\int_{-\infty}^{x}\hat{\eta}(\rho(y,\tau),(\rho
u)(y,\tau))\ dy\bigg)\ dx\bigg|.
\end{aligned}
\end{equation}

Applying Lemma $\ref{lem3.1}$, it is easy to get
\begin{eqnarray}\label{3.44}
\int_{0}^{t}\int_{K}|\eta^{\star}(\rho,\rho(u-\bar{u}))|\ dxd\tau\leq
C(t).
\end{eqnarray}

Now Cauchy-Schwartz inequality and  \eqref{3.55} lead to
\begin{eqnarray}\label{3.38}
\varepsilon\int_{0}^{t}\int_{K}\rho^{\alpha+\theta}|u_x|\ dxd\tau\
&\leq& C\varepsilon\int_{0}^{t}\int_{K}\rho^{\alpha}u_x^2\ dxd\tau\
+C\varepsilon\int_{0}^{t}\int_{K}\rho^{\alpha+2\theta}\ dxd\tau\nonumber\\
&\leq& C(t)+C(t)\int_{0}^{t}\int_{K}\rho^{\frac{\gamma-1}{2}}\ dxd\tau\nonumber\\
&\leq& C(t).
\end{eqnarray}

Noticing that $3\beta>3\alpha-2$ and \eqref{3.55}, we have
\begin{eqnarray}\label{3.39}
\varepsilon\int_{0}^{t}\int_{K}\rho^{\alpha}|u_x||u|\ dxd\tau
&\leq&\frac12\varepsilon\int_{0}^{t}\int_{K}\rho^{\alpha}u_x^2\ dxd\tau\
+\frac12\varepsilon\int_{0}^{t}\int_{K}\rho^{\alpha}u^2\ dxd\tau\nonumber\\
&\leq&
C(t)+C(\delta)\varepsilon^3\int_{0}^{t}\int_{K}\rho^{3\alpha-2}\ dxd\tau
+\delta\int_{0}^{t}\int_{K}\rho|u|^3\ dxd\tau\nonumber\\
&\leq& C(t)+\varepsilon^3\int_{0}^{t}\int_{K}\rho^{3\beta}\ dxd\tau
+\delta\int_{0}^{t}\int_{K}\rho|u|^3\ dxd\tau\nonumber\\
&\leq& C(t) +\delta\int_{0}^{t}\int_{K}\rho|u|^3\ dxd\tau
\end{eqnarray}
where $\delta$ is small enough which will be determined later.

Now we estimate the last term on the right hand side
of \eqref{3.37}. \eqref{1.1} implies that 
\begin{eqnarray}\label{3.46}
(\rho u-\rho u^-)_t+(\rho u^2+p(\rho)-\rho
uu^-)_x=\varepsilon(\rho^{\alpha}u_x)_x.
\end{eqnarray}
Integrating \eqref{3.46} over $[0,t]\times(-\infty,x)$ for $x\in K$,
we obtain
\begin{eqnarray}\label{3.47}
\int_{-\infty}^{x}\rho( u- u^-)\ dy&=&\int_{-\infty}^{x}\rho_0( u_0-
u^-)\ dy-\int_{0}^{t}(\rho u^2+p(\rho)-\rho
uu^--p(\rho^-))\ d\tau\nonumber\\
&+&\varepsilon\int_{0}^{t}\rho^{\alpha}u_x\ d\tau.
\end{eqnarray}
On the other hand,
\begin{equation}\label{3.45}
\begin{aligned}
&\displaystyle\bigg|\int_{-\infty}^{x}\hat{\eta}(\rho(y,\tau),(\rho
u)(y,\tau))\ dy\bigg|\\
&\displaystyle\leq \bigg|\int_{-\infty}^{x}(\hat{\eta}(\rho,\rho
u)-\beta\rho^{\theta+1}(u-\bar{u}))\ dy\bigg|
+\bigg|\int_{-\infty}^{x}\beta\rho^{\theta+1}(u-u^-)\ dy\bigg|\\
&\displaystyle\leq \bigg|\int_{-\infty}^{x}r(\rho,\rho(u-u^-))\ dy\bigg|
+\bigg|\int_{-\infty}^{x}\beta(\rho^{\theta}-(\rho^-)^{\theta})\rho(u-u^-)\ dy\bigg|\\
&\displaystyle\qquad+\beta(\rho^-)^{\theta}\bigg|\int_{-\infty}^{x}\rho(u-u^-)\ dy\bigg|\\
&\displaystyle\leq
C(t)+\beta(\rho^-)^{\theta}\bigg|\int_{-\infty}^{x}\rho(u-u^-)\ dy\bigg|,
\end{aligned}
\end{equation}
which, together with \eqref{3.51}, Lemma$\ref{lem3.1}$-Lemma$\ref{lem3.3}$ and
\eqref{3.47}, implies that 
\begin{eqnarray}\label{3.48}
\int_{K}\bigg|\int_{-\infty}^{x}\hat{\eta}(\rho(y,\tau), m(y,\tau))dy\bigg|\ dx\leq C(t).
\end{eqnarray}

Now, if one chooses $\delta$ small enough, then substitutes \eqref{3.48}, \eqref{3.44}, \eqref{3.39} and \eqref{3.38}
into \eqref{3.37}, the proof of Lemma 2.4 follows.

\begin{remark}
In the uniform estimates above, we have required that
$\frac23\leq\alpha\leq\gamma$, and the initial functions
$\left({\rho^\varepsilon_0(x),u^\varepsilon_0(x)}\right)$ satisfy

(i) $\displaystyle\rho^\varepsilon_0(x)>0,\quad
\int_{\mathbb{R}}\rho^\varepsilon_0(x)|u^\varepsilon_0(x)-\bar{u}(x)|dx\leq
M_0<\infty$;

(ii) The total mechanical energy with respect to
$\left({\bar\rho,\bar{u}}\right)$ is finite:
$$\int_{\mathbb{R}}\frac{1}{2}\rho^\varepsilon_0\left({x}\right)\left|{u^\varepsilon_0(x)-\bar{u}(x)}\right|^2+e^{\ast}\left({\rho^\varepsilon_0(x),\bar\rho(x)}\right)dx=:E_0<\infty;$$

(iii)
$\displaystyle\varepsilon^2\int_{\mathbb{R}}\frac{\left|{\rho^\varepsilon_{0x}\left({x}\right)}\right|^2}{\rho^\varepsilon_{0}\left({x}\right)^{3-2\alpha}}dx\leq
E_1<\infty$.

(iv) $M_0$, $E_0$, $E_1$ are independent of $\varepsilon$.
\
These conditions are essential parts of Condition 1 in section 1. We remark here that the limit of the functions satisfying the
conditions i)-iv) is very general, including a wide class of $L^\infty$ functions with 
finite energy and may contain vacuum. It is obvious that the above limit can serve as the 
initial data of
isentropic gas dynamics for the existence of finite energy solutions. We also note that the condition
(iii) is slightly weaker than the corresponding one in \cite{Chen6} near vacuum. 
We refer to \cite{Chen6} for further details.
\end{remark}

%%%%%%%%%%%%%%%%%%%%%%%%%%%%%%%%%%%%%%%%%%%%%%%%%%%%%%%%%%%%%%%%%%%%%%%%%%%%%%%%%%%%%%%%%%%%%%%%%%%%%%%
\section{$H_{loc}^{-1}(\mathbb{R}_+^2)-$ Compactness}

In this section we will use the uniform estimates obtained
in the previous section to prove the following key Lemma, which
states the
$H_{loc}^{-1}(\mathbb{R}_+^2)-$compactness of the approximate 
solution sequence.

\begin{lemma}\label{lem4.1}
Let $\frac23\leq\alpha \leq \gamma$, $\psi\in C_0^2(\mathbb{R})$,
$(\eta^{\psi},q^{\psi})$ be a weak entropy pair generated by $\psi$.
Then for the solutions $(\rho^{\varepsilon},u^{\varepsilon})$ with
$m^{\varepsilon}=\rho^{\varepsilon}u^{\varepsilon}$ of Navier-Stokes
equations (1.1)--(1.2), the following sequence
\begin{eqnarray}\label{4.1}
\eta^{\psi}(\rho^{\varepsilon},m^{\varepsilon})_t+q^{\psi}(\rho^{\varepsilon},m^{\varepsilon})_x
\ \mbox{are compact in $H_{loc}^{-1}(\mathbb{R}_+^2)$}
\end{eqnarray}
\end{lemma}

\noindent\textbf{Proof}. In order to prove this lemma, we first
introduce the following results for the entropy pair
$(\eta^{\psi},q^{\psi})$ generated by $\psi\in C_0^2(\mathbb{R})$,
and see \cite{Chen6} for details.

For a $C^2$ function $\psi:\mathbb{R}\rightarrow\mathbb{R}$,
compactly supported on the interval $[a,b]$, we have
\begin{eqnarray}\label{4.10}
supp\{\eta^{\psi}\},supp \{q^{\psi}\}\subset \left\{{(\rho,m)=(\rho,\rho
u):u+\rho^{\theta}\geq a,\quad u-\rho^{\theta}\leq b}\right\}.
\end{eqnarray}
Furthermore, there exists a constant $C_{\psi}>0$ such that, for any
$\rho\geq0$ and $u\in\mathbb{R}$, we have

(i) For $\gamma\in(1,3]$,
\begin{eqnarray}\label{4.11}
|\eta^{\psi}(\rho,m)|+|q^{\psi}(\rho,m)|\leq C_{\psi}\rho.
\end{eqnarray}

(ii) For $\gamma\in(3,+\infty)$,
\begin{eqnarray}\label{4.12}
|\eta^{\psi}(\rho,m)|\leq C_{\psi}\rho,\ \ |q^{\psi}(\rho,m)|\leq
C_{\psi}(\rho+\rho^{\theta+1}).
\end{eqnarray}

(iii) If $\eta^{\psi}$ is considered as a function of $(\rho,m)$,
$m=\rho u$, then
\begin{eqnarray}\label{4.13}
|\eta^{\psi}_m(\rho,m)|+|\rho\eta^{\psi}_{mm}(\rho,m)|\leq C_{\psi},
\end{eqnarray}
 and, if $\eta^{\psi}_m$ is considered as a
function of $(\rho, u)$, then
\begin{eqnarray}\label{4.14}
|\eta^{\psi}_{mu}(\rho,\rho
u)|+|\rho^{1-\theta}\eta^{\psi}_{m\rho}(\rho,\rho u)|\leq C_{\psi}.
\end{eqnarray}

Now we are going to prove the lemma. 

A direct computation on $\eqref{1.1}_1\times\eta^{\psi}_{\rho}(\rho^{\varepsilon},m^{\varepsilon})
+\eqref{1.1}_2\times\eta^{\psi}_m(\rho^{\varepsilon},m^{\varepsilon})$ 
gives 
\begin{eqnarray}\label{4.2}
&\displaystyle\eta^{\psi}(\rho^{\varepsilon},m^{\varepsilon})_t+q^{\psi}(\rho^{\varepsilon},m^{\varepsilon})_x
=\varepsilon(\eta^{\psi}_m(\rho^{\varepsilon},m^{\varepsilon})(\rho^{\varepsilon})^{\alpha}u_x^{\varepsilon})_x
-\varepsilon\eta^{\psi}_{mu}(\rho^{\varepsilon},m^{\varepsilon})(\rho^{\varepsilon})^{\alpha}(u_x^{\varepsilon})^2\nonumber\\
&\displaystyle\qquad\qquad\qquad\quad-\varepsilon\eta^{\psi}_{m\rho}(\rho^{\varepsilon},m^{\varepsilon})(\rho^{\varepsilon})^{\alpha}\rho_x^{\varepsilon}u_x^{\varepsilon}
\end{eqnarray}

Let $K\subset\mathbb{R}$ be compact, using \eqref{4.14} and Cauchy-Schwartz
inequality, we have
\begin{equation}\label{4.4}
\begin{aligned}
&\displaystyle\varepsilon\int_{0}^{t}\int_{K}|\eta^{\psi}_{mu}(\rho^{\varepsilon},m^{\varepsilon})(\rho^{\varepsilon})^{\alpha}
|(u_x^{\varepsilon})^2
+|\eta^{\psi}_{m\rho}(\rho^{\varepsilon},m^{\varepsilon})(\rho^{\varepsilon})^{\alpha}\rho_x^{\varepsilon}u_x^{\varepsilon}|\ dxdt\\
&\displaystyle\leq
C\varepsilon\int_{0}^{t}\int_{K}(\rho^{\varepsilon})^{\alpha}(u_x^{\varepsilon})^2\ dxd\tau
+C\varepsilon\int_{0}^{t}\int_{K}(\rho^{\varepsilon})^{\alpha+\gamma-3}(\rho_x^{\varepsilon})^2\ dxd\tau\\
&\leq C(t).
\end{aligned}
\end{equation}
This implies that
\begin{eqnarray}\label{4.5}
-\varepsilon\eta^{\psi}_{mu}(\rho^{\varepsilon},m^{\varepsilon})(\rho^{\varepsilon})^{\alpha}(u_x^{\varepsilon})^2
-\varepsilon\eta^{\psi}_{m\rho}(\rho^{\varepsilon},m^{\varepsilon})(\rho^{\varepsilon})^{\alpha}\rho_x^{\varepsilon}u_x^{\varepsilon}
\ \mbox{are bounded in $L^1([0,T]\times K)$},
\end{eqnarray}
and thus it is compact in
$W_{loc}^{-1,p_1}(\mathbb{R}_+^2)$, for $1<p_1<2$.

Moreover, noticing
$|\eta^{\psi}_{m}(\rho^{\varepsilon},\rho^{\varepsilon}u^{\varepsilon})|\leq
C_{\psi}$, we have
\begin{equation}\label{4.6}
\begin{aligned}
&\displaystyle\int_{0}^{t}\int_{K}\bigg(\varepsilon
\eta^{\psi}_m(\rho^{\varepsilon},m^{\varepsilon})
(\rho^{\varepsilon})^{\alpha}u_x^{\varepsilon}\bigg)^{\frac{4}{3}}\ dxdt\\
&\displaystyle\leq\int_{0}^{t}\int_{K}\varepsilon^{\frac{4}{3}}
(\rho^{\varepsilon})^{^{\frac{4\alpha}{3}}}|u_x^{\varepsilon}|^{\frac{4}{3}}dxdt\\
&\displaystyle\leq C\varepsilon^{\frac{4}{3}}\int_{0}^{t}\int_{K}
(\rho^{\varepsilon})^{\alpha}|u_x^{\varepsilon}|^2dxdt+C\varepsilon^{\frac{4}{3}}\int_{0}^{t}\int_{K}
(\rho^{\varepsilon})^{2\alpha}dxdt\\
 &\displaystyle\leq C(T,K)\varepsilon^{\frac{1}{3}}+C\varepsilon^{\frac{4}{3}}\int_{0}^{t}\int_{K}
(\rho^{\varepsilon})^{\gamma+1}dxdt\\
&\displaystyle\leq C(T,K)\varepsilon^{\frac{1}{3}}\rightarrow  0\ \ \mbox{as}\
\varepsilon\rightarrow 0.
\end{aligned}
\end{equation}
Then\eqref{4.6} and \eqref{4.5} yield that
\begin{eqnarray}\label{4.7}
\eta^{\psi}(\rho^{\varepsilon},m^{\varepsilon})_t+q^{\psi}(\rho^{\varepsilon},m^{\varepsilon})_x
\ \mbox{are compact in $W_{loc}^{-1,p_2}(\mathbb{R}_+^2)$} \ \
\mbox{for some}\ 1<p_2<2.
\end{eqnarray}

On the other hand, using the estimates in \eqref{4.11}--\eqref{4.12}
and Lemma $\ref{lem3.1}$--Lemma $\ref{lem3.4}$, we have
\begin{eqnarray}\label{4.8}
\eta^{\psi}(\rho^{\varepsilon},m^{\varepsilon}), q^{\psi}(\rho^{\varepsilon},m^{\varepsilon})
\ \mbox{are uniformly bounded in $L_{loc}^{p_3}(\mathbb{R}_+^2)$} \ \
\mbox{for}\ p_3>2,
\end{eqnarray}
where $p_3=\gamma+1>2$ when $\gamma\in(1,3]$; and
$p_3=\frac{\gamma+\theta}{1+\theta}>2$ when $\gamma>3$. This yields
that,
\begin{eqnarray}\label{4.9}
\eta^{\psi}(\rho^{\varepsilon},m^{\varepsilon})_t+q^{\psi}(\rho^{\varepsilon},m^{\varepsilon})_x
\ \mbox{are uniformly bounded in $W_{loc}^{-1, p_3}(\mathbb{R}_+^2)$}.
\end{eqnarray}
Then \eqref{4.7} and \eqref{4.9} implies Lemma $\ref{lem4.1}$.

%%%%%%%%%%%%%%%%%%%%%%%%%%%%%%%%%%%%%%%%%%%%%%%%%%%%%%%%%%%%%%%%%%%%%%%%%%%%%%%%%%%%%%%%%%%%%%%%%%%%%%%

\section{Proof of Theorem}

\noindent\textbf{Proof of Theorem \ref{th1.1}}. From  Lemmas
\ref{lem3.1}--\ref{lem3.4} and  the compactness estimate Lemma
\ref{lem4.1}, we have verified the conditions (i)-(iii) of  Theorem
\ref{th4.1}  for the sequence of solutions
$(\rho^\varepsilon, m^\varepsilon)$. Basing on Theorem
\ref{th4.1},  there is a subsequence
$(\rho^{\varepsilon},m^{\varepsilon})$(still denoted as $(\rho^{\varepsilon},m^{\varepsilon})$) and a pair of
measurable functions $(\rho, m)$ such that
\begin{eqnarray}\label{4.1-2}
(\rho^{\varepsilon},m^{\varepsilon})\rightarrow (\rho, m),\ \ a.e \
\varepsilon\rightarrow0.
\end{eqnarray}

It is easy to check that $(\rho,m)$ is a finite-energy entropy
solution $\left({\rho,m}\right)$ to the Cauchy problem (\ref{1.3})
with initial data $\left({\rho_0,\rho_0 u_0}\right)$ for the
isentropic Euler equations with $\gamma>1$. Therefore, the proof of Theorem
\ref{th1.1} is completed.

\

\noindent {\bf Acknowledgments} The research of FMH was supported in
part by NSFC Grant No. 10825102 for Outstanding Young scholars,
National Basic Research Program of China (973 Program),
No.2011CB808002. The research of RHP was supported by National
Science Foundation under grant DMS-0807406.

\

\end{document}